\documentclass[a4paper]{article}
\usepackage{times}
\usepackage[T1]{fontenc}
\usepackage[left=1.5cm,right=1.5cm,top=1.5cm,bottom=1.7cm]{geometry}
\usepackage[utf8]{inputenc}
\usepackage{algorithm}
\usepackage[noend]{algpseudocode}
\usepackage{amsmath,amssymb}
\usepackage{subcaption}
\usepackage{hyperref}
\usepackage{pgfplots}
\usepackage{tikz}
\usepackage{graphicx}
\usepackage[symbol]{footmisc}

\let\oldbibliography\thebibliography
\renewcommand{\thebibliography}[1]{%
  \oldbibliography{#1}%
  \setlength{\itemsep}{-0.1cm}%
}

\begin{document}
\title{Multigrid-reduction-in-time for Eddy Current problems}

\author{\textbf{Stephanie Friedhoff}$^{1}$, \textbf{Jens Hahne}$^{1,*}$ and \textbf{Sebastian Schöps}$^{2}$}

\date{\small
$^1$ Department of Mathematics, Bergische Universität Wuppertal, Gaußstr. 20, 42119 Wuppertal, Germany\\
$^2$ Centre for Computational Engineering, Technische Universität Darmstadt, Dolivostr. 15, 64293 Darmstadt, Germany
}

\maketitle
\footnotetext[1]{Corresponding author: e-mail \texttt{jens.hahne@math.uni-wuppertal.de}}

\section*{Abstract}
\noindent Parallel-in-time methods have shown success for reducing the simulation time of many time-dependent problems.
Here, we consider applying the multigrid-reduction-in-time (MGRIT) algorithm to a voltage-driven eddy current model problem.

\section{Introduction}
The simulation of electrical machines, such as synchronous and induction machines, transformers or cables, is an established procedure in industry to improve the product design, e.\,g., to prevent eddy current losses.
In this low-frequency regime, the eddy current problem is typically used. It is an approximation of Maxwell's equations for devices where the displacement current can be neglected with respect to the source currents \cite{Schmidt_2008aa}. For such devices, the model governs the evolution of electromagnetic fields and is, for a voltage-driven system, coupled with an additional equation, resulting in the following system for unknown magnetic vector potential $\vec{A}: \Omega \times \mathcal{I} \rightarrow \mathbb{R}^3$ and the electric current $i_s: \mathcal{I} \rightarrow \mathbb{R}$:
\begin{subequations}\label{eq:system}
\begin{align}
\sigma \partial_t \vec{A} + \nabla \times \big(\nu(\|\nabla\times\vec{A}\|) \nabla\times\vec{A}\big) - \vec{\chi}_{s} i_s &= 0, \label{eq:ecp} \\
\frac{d}{dt}\int_{\Omega} \vec{\chi}_{s}\cdot\vec{A}\; dV & = \upsilon_{s},\label{eq:vi}
\end{align}
\end{subequations}
with homogeneous Dirichlet condition $\vec{A}\times\vec{n}=0$ with normal vector $\vec{n}$ on $\partial\Omega$, initial value $\vec{A}|_{t_0}\equiv0$, and where $ \Omega $ denotes the spatial domain (here, the tube region depicted in Fig.~\ref{fig:model}) and $ \mathcal{I} = (t_0, t_{\text{end}}]$ is the time interval. The electrical conductivity $ \sigma \geq 0$ is only non-zero in the tube region $\Omega_2 \; (10$ MS/m$)$, and the magnetic reluctivity $ \nu $ is modeled by a monotone cubic spline curve in $\Omega_2$ and by vacuum $(1 / \mu_0)$ in $\Omega_0$ and $\Omega_1$. The winding function $\vec {\chi}_s: \Omega \rightarrow \mathbb{R}^3$ represents a stranded conductor in the model \cite{Schoeps_2013}. Equation \eqref{eq:vi} establishes a relationship between the so-called flux linkage, i.\,e., the spatially integrated time derivative of the magnetic vector potential and the pulsed voltage $v_s(t)= 0.25\;p(t)$ V, which is given by
\begin{equation}\label{eq:input_pwm}
	p(t)= \begin{cases}
	\mathrm{sign}\left[\sin\left(\dfrac{2\pi}{T} t\right)\right],\ & s_n(t)-\left|\sin\left(\dfrac{2\pi}{T}t\right)\right|<0,\\
	0,\ & \mathrm{otherwise,}
	\end{cases}
\end{equation}
where $s_n(t)=n/Tt -\left\lfloor n/Tt\right\rfloor$ is the common sawtooth pattern, with $n=200$ teeth and period $T=0.02\;$s \cite{Kulchytska-Ruchka_2018ac}. Discretizing in space using edge shape functions yields a system of index-$1$ differential-algebraic equations (DAEs). Please note, that in this particular case, the symmetry in the $z$-direction is exploited and the problem is solved only on a 2D section in the $x$-$y$ plane. This semi-discrete system can be integrated using the backward Euler method, resulting for each time step $t_j$ in a nonlinear system of the form $\Phi(\mathbf{u}_j)=\mathbf{g}_j$, with $\mathbf{u}^{\!\top}_j=(\mathbf{a}^{\!\top}, i)$ where $\mathbf{a}$ is the vector of discrete vector potentials and $i$ is an approximation of the current. Collecting all $\mathbf{u}_j$'s and $\mathbf{g}_j$'s in vectors $\mathbf{u}$ and $\mathbf{g}$, respectively, we obtain the space-time system $A(\textbf{u}) = \textbf{g}$, where each block row corresponds to one time step. 

\begin{figure}[h!]
\begin{minipage}[t][][c]{0.48\textwidth}
	\begin{tikzpicture}[scale=.7]
		\draw[semithick,fill=gray!70] (0,0) ellipse (0.2 and 0.46);
		\draw[semithick,fill=white] (0,0) ellipse (0.1 and 0.3);
		\draw[line width=2pt,color=black!80] (-1.3,0) -- (.1,0);
		\draw[line width=2pt,color=black!80] (.5,0) -- (4.3,0);
		\draw[semithick] (0,0.46) -- (3,0.46);
		\draw[semithick] (0,-0.46) -- (3,-0.46);
		\draw[semithick] (3,-0.46) arc(-90:90:0.2 and 0.46);
		\filldraw[fill=gray!70, draw=black]
			(0,-.46) -- (3,-.46) arc (-90:90:0.2 and 0.46) -- (0,.46) arc (90:-90:0.2 and 0.46);
			
		\draw[semithick,->,>=latex] (-.5,-1.3) -- (.06,-1.65) node[right=-1pt] {$\small x$};
		\draw[semithick,->,>=latex] (-.5,-1.3) -- (-.5,-.6) node[above right=-3pt] {$\small y$};
		\draw[semithick,->,>=latex] (-.5,-1.3) -- (-1.2,-1.3) node[left=-1pt] {$z$};
		
		\draw[semithick,fill=gray!70] (0+7,0) circle (40pt);
		\draw[semithick,fill=white] (0+7,0) circle (25pt);
		\draw[semithick,fill=black!80] (0+7,0) circle (10pt);
		
		\draw[semithick,->,>=latex] (-2+7,-1.3) -- (-1.4+7,-1.3) node[right=-2pt] {$x$};
		\draw[semithick,->,>=latex] (-2+7,-1.3) -- (-2+7,-.6) node[above right=-3pt] {$y$};
	\end{tikzpicture}
\caption{\label{fig:model} Coaxial cable model and its cross section. The inner, dark grey region $\Omega_0$ models the wire, the white region $\Omega_1$ an insulator and the outer, light grey region $\Omega_2$ the conducting shield \cite{FEMM}.}
\end{minipage}
\hfill
	\begin{minipage}[t][][c]{0.48\textwidth}
\begin{algorithm}[H] \small
\caption{\label{alg:mgrit}MGRIT($A, \mathbf{u}, \mathbf{g}$)}
\begin{algorithmic}[1]
\State Apply $FCF$-relaxation to $A_1(\textbf{u}^{(1)})= \textbf{g}^{(1)}$
\State Inject the approximation and its residual to the coarse grid: 
\Statex \qquad $\textbf{u}^{(2)} = R_I(\textbf{u}^{(1)}), $
\Statex \qquad $\textbf{g}^{(2)} = R_{I}(\textbf{g}^{(1)}-A_1\textbf{u}^{(1)})$
\State Solve $A_{2}(\textbf{v}^{(2)})= A_{2}(\textbf{u}^{(2)}) + \textbf{g}^{(2)}$
\State Compute the error approximation: $\mathbf{e} = \mathbf{v}^{(2)} - \mathbf{u}^{(2)}$
\State Correct using ideal interpolation: $\textbf{u}^{(1)} = \textbf{u}^{(1)} + P(\mathbf{e}) $
\end{algorithmic}
\end{algorithm}
\end{minipage}
\end{figure}


\section{Multigrid reduction in time (MGRIT)}
The multigrid-reduction-in-time (MGRIT) algorithm \cite{Falgout_etal_2014} is a highly parallel, iterative method for solving time-dependent problems. Instead of sequential time stepping, MGRIT is a parallel-in-time method that allows for solving multiple time steps at once by using a multilevel hierarchy of temporal grids. For a given time grid, define a splitting of all time points into $C$- and $F$-points, such that every $m$-th point is a $C$-point, defining the next coarser grid, and all others are $F$-points. MGRIT uses block relaxation that alternates between $F$- and $C$-points. Relaxation on $F$-points, called $F$-relaxation, propagates the solution from each $C$-point to all $F$-points up to the next $C$-point. Similarly, $C$-relaxation updates the solution at all $C$-points. Further, we define two grid transfer operators, injection as the restriction method and ``ideal interpolation'' is injection followed by an $F$-relaxation. The resulting two-level MGRIT algorithm is given in Algorithm \ref{alg:mgrit}, where $A_l(\mathbf{u}^{(l)}) = \mathbf{g}^{(l)}$ denote the nonlinear space-time systems of equations on levels $l = 1, 2$ with each block row corresponding to one time step; multilevel schemes, e.\,g., $V$- and $F$-cycles \cite{UTrottenberg_etal_2001} can be defined by applying the two-level method recursively to the system in Step $3$. 


\section{Numerical results}
We apply the MGRIT algorithm to the coupled system \eqref{eq:system} on the space-time domain $\Omega \times (0,0.04]$ with $\Omega = \Omega_0 \cup \Omega_1 \cup \Omega_2$ (Fig.~\ref{fig:model}). The problem is discretized using linear edge shape functions with $2269$ degrees of freedom in space and on a fine equidistant time grid with $2^{14}$ intervals using backward Euler to resolve the pulses. We consider MGRIT $V$- and $F$-cycle variants with different numbers of grid levels $(L = 3, 4, 5)$ and various coarsening factors $(m = 64, 16, 8)$, each with a coarsest grid consisting of four time points. The system on the coarsest grid is solved with time stepping and all nonlinear spatial solves are performed using Newton's method.
For numerical results, an Intel Xeon Phi cluster with $272$ $1.4$ GHz Intel Xeon Xhi processors is used. The MGRIT algorithm was implemented in Python using the Message Passing Interface (MPI). 

\begin{figure}[h!]
	\begin{minipage}[c][14\baselineskip][c]{0.45\textwidth}
		\begin{center}
				\vfill
			\begin{tabular}{@{}l|ccc@{}}
				\hline
					& $L=3$ & $L=4$ & $L=5$ \\
				\hline
					$V$-cycle $FCF$ & $7$ & 9 & 9 \\
					$F$-cycle $FCF$ & 7 & 9 & 8 \\
				\hline
			\end{tabular}
					\vfill
		\end{center}
		\caption{\label{tab:mgrit_iter}Iteration counts of MGRIT variants for solving the coupled system \eqref{eq:system} to a space-time residual norm $\|r\| < 10^{-6}$.}
	\end{minipage}
	\hfill
	\begin{minipage}[c][14\baselineskip][c]{0.48\textwidth}
		\includegraphics[width=\textwidth]{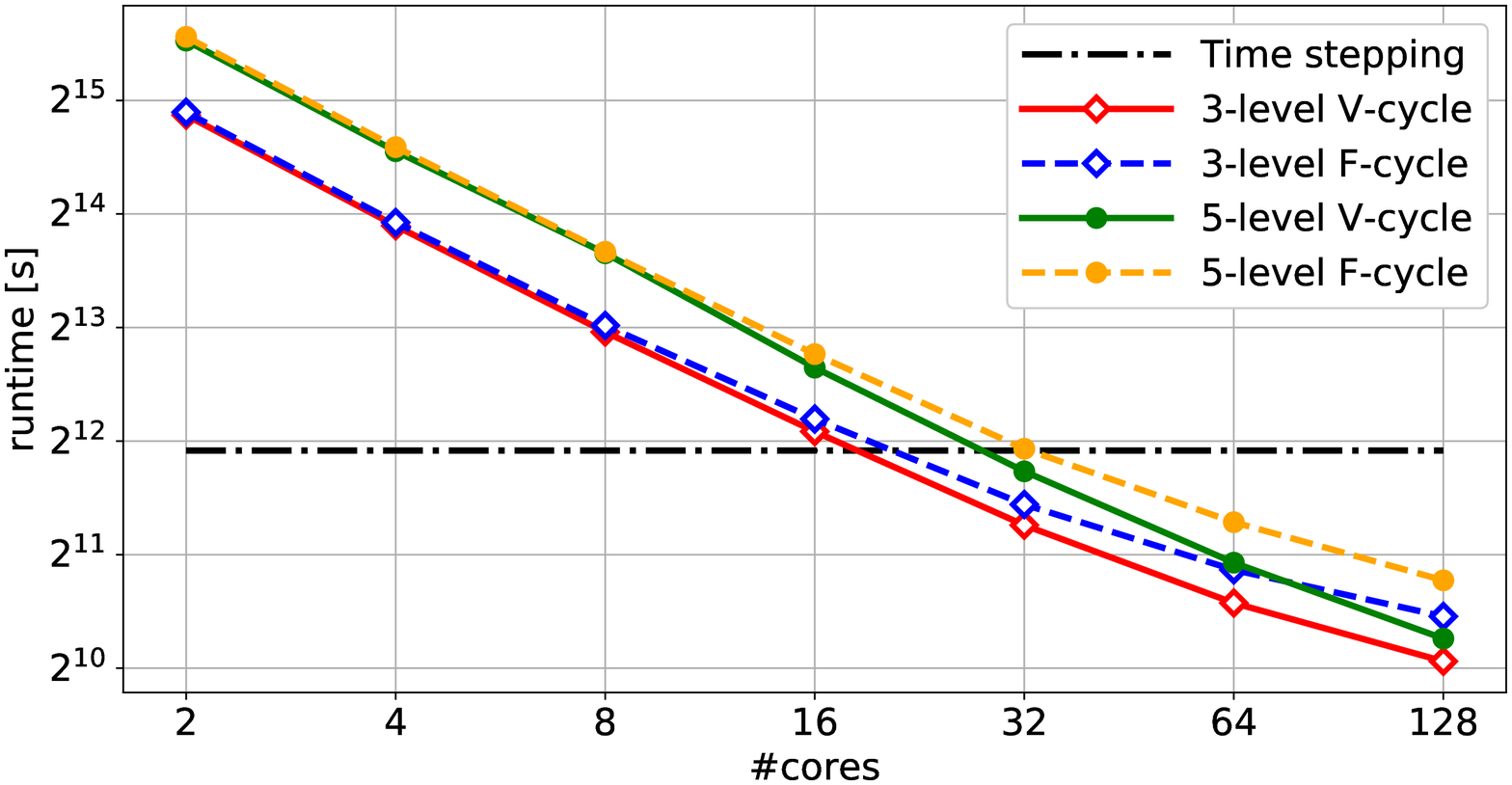}
		\caption{\label{fig:strong_scaling}Time stepping on one processor and strong scaling results for MGRIT variants applied to the coupled system. }
	\end{minipage}
\end{figure}

The results in Tab. \ref{tab:mgrit_iter} show that all six MGRIT variants have similar convergence behavior. While iteration counts of the two cycle types are the same for three and four levels, for five levels, the number of iterations of the $F$-cycle MGRIT algorithm is smaller than the $V$-cycle variant. Note, however, that the cost of one $F$-cycle iteration grows with increasing numbers of grid levels. Fig. \ref{fig:strong_scaling} shows strong scaling results for the three- and five-level variants as well as the runtime of time stepping on one processor for reference purposes. We see good scaling behavior for all four schemes, with a speedup of up to a factor of about 4.2 over sequential time stepping. On smaller numbers of processors, runtimes of $F$- and $V$-cycles are about the same, while for larger processor counts $V$-cycles are faster than $F$-cycles, due to the higher communication requirements of $F$-cycles.

\paragraph{Acknowledgement} The work is supported by the Excellence Initiative of the German Federal and State Governments, the Graduate School of Computational Engineering at TU Darmstadt, and the BMBF (project PASIROM; grants 05M18RDA and 05M18PXB).


\end{document}